\documentclass[11pt]{article}
\usepackage{amssymb,latexsym}
\def\RR{\mathbb{R}}
\def\ZZ{\mathbb{Z}}

\def\NN{\mathbb{N}}
\def\PP{\mathbb{P}}
\def\VV{\mathbb{V}}
\begin{document}

\title{\bf Approximation by Baskakov quasi-interpolants}
\author{Paul Sablonni\`ere, INSA \& IRMAR, Rennes}
\maketitle

\begin{abstract}
Baskakov operators and their inverses can be expressed as linear differential operators on polynomials. Recurrence relations are given for the computation of these coefficients. They allow the construction of  the associated Baskakov  quasi-interpolants (abbr. QIs).
Then asymptotic results are provided for the determination of the convergence orders of these new quasi-interpolants. Finally some results on the computation of these QIs and the numerical approximation of functions defined  on the positive real half-line are  illustrated by some numerical examples.
\end{abstract}
%-----------------------------------------------------------------------------------------------------------------------------------
%**************************************1./INTRODUCTION/*******************************************
%-----------------------------------------------------------------------------------------------------------------------------------
\section{Introduction}
In the present paper, the general method developed by various authors for the construction of  quasi-interpolants (abbr. QIs) of Bernstein and other types, is applied to Baskakov operators \cite{B} defined by
$$
\mathcal{V}_nf(x):=\sum_{k\ge 0} f_k \, v_{k,n}(x),\quad f_k:=f\left(\frac kn\right)
$$
where
$$
v_{k,n}(x):= {n+k-1 \choose k} x^k (1+x)^{-(n+k)},\, k\ge 0
$$
 We thus complete the results given by P. Mache and M.W. M\"uller in \cite{MM2}. The method also works for Durrmeyer type Baskakov operators and will be developed elsewhere. \\
This method can be summarized as follows : let $\{Q_n, n\in \NN\}$ be a sequence of linear  operators defined on some functional space $\mathcal{F}$ with values in a finite-dimensional subspace $\mathcal{P}_n$ of algebraic (or trigonometric) polynomials. Assuming that for all $n\in \NN$, $Q_n$ is an isomorphism of 
$\mathcal{P}_n$ preserving the degree, i.e. $Q_np\in \mathcal{P}_m$ for any $p\in \mathcal{P}_m, 0\le m\le n$, very often it admits a representation in that space as a differential operator of the form
$$
Q_n=\sum_{r=0}^n \beta_r^{(n)} \mathcal{D}^r
$$
where $\mathcal{D}$ is a (simple) linear differential operator, and $\beta_r^{(n)}(x)$ is a polynomial of degree at most $r$. In most cases, the inverse $P_n:=Q_n^{-1}$ of $Q_n$ has also a representation of the same form
$$
P_n=\sum_{r=0}^n \alpha_r^{(n)} \mathcal{D}^r.
$$
In general, both families of polynomial coefficients satisfy a recurrence relation. This has been proved 
\cite{Sab3} for Bernstein and Sz\'asz-Mirakyan  operators  and their associated Kantorovich  and Durrmeyer versions, and also for Weierstrass operators \cite{Sab6}.\\

Introducing the truncated inverse of order $0\le r\le n$ of $P_n$
$$
P_n^{(r)}=\sum_{k=0}^r \alpha_k^{(n)} \mathcal{D}^k
$$
one can associate with $Q_n$ either the left quasi-interpolants $\{Q_n^{(r)}, 0\le r\le n\}$ defined by
$$
Q_n^{(r)}p:=P_n^{(r)} Q_n p=\sum_{k=0}^r \alpha_k^{(n)}(x) \mathcal{D}^k Q_n p,\quad \forall p\in \mathcal{P}_n
$$
or the right quasi-interpolants $\{Q_n^{[r]}, 0\le r\le n\}$ defined by
$$
Q_n^{[r]}p:=Q_nP_n^{(r)} p=\sum_{k=0}^r Q_n(\alpha_k^{(n)}\mathcal{D}^k p),\quad \forall p\in \mathcal{P}_n
$$
By construction, these operators are exact on $\mathcal{P}_r$ and, in general, they can be extended to the functional space $\mathcal{F}$. In virtue of classical theorems in approximation theory, this procedure greatly improves the convergence order of the initial operator $Q_n$ in the space $\mathcal{F}$.

After its introduction in \cite{Sab1,Sab2} for Bernstein operators (see also \cite{MM1}), it has been extended by several authors to various classical linear approximation operators. 
For example  by P. Mache and M.W. M\"uller \cite{MM2} to Baskakov operators, by A.T. Diallo \cite{D1,D3,D4} and the author \cite{Sab3} to Sz\'asz-Mirakyan operators, and more recently \cite{Sab6} to Weierstrass operators.
This general method also works for multivariate extensions of these operators (see e.g. \cite{Sab5} for Bernstein operators on triangles).\\

Here is a brief outline of the paper. In Section 2, we construct the polynomial coefficients of the differential operators representing the Baskakov operator and its inverse on the space of polynomials. Theorem 1 and 2 give the recurrence relations for the computation of these polynomials and their asymptotic behaviour when $n\to \infty$. In Section 3, some additional results are given on the norms of Baskakov QIs $\mathcal{V}_n^{(r)}$ and a Voronovskaya type theorem about their asymptotic convergence orders on smooth functions. Section 4 describes some practical methods for the effective computation of these operators.  Section 5 presents some numerical examples of approximation of functions by Baskakov quasi-interpolants. Finally, in Section 6, are listed the polynomial coefficients that are needed for the construction of the QIs $\mathcal{V}_n^{(r)}$ for $5\le r\le 11$.  \\
We want to emphasize the fact that {\em the function to be approximated is only known by its values on uniform partitions} of step $1/n$ of the positive real axis $\RR_+=[0,+\infty)$. Therefore the Baskakov QIs cannot be compared with operators using specific points like e.g. zeros of orthogonal polynomials or similar systems of abscissas for the evaluation of the function. In a forthcoming paper \cite{Sab7}, we will study the quadrature formulas on $\RR_+$ generated by integration of those Baskakov QIs.\\

{\bf Notations.} The rising factorials are written
$$
(n)_r=\frac{(n+r-1)!}{(n-1)!}=n(n+1)\ldots(n+r-1)
$$
The (positive) Stirling numbers of the first kind (see e.g. \cite{A}, chapter 1) are defined by
$$
(x)_n=\sum_{k=0}^n (-1)^{n-k} s(n,k) x^k
$$
The falling factorials are written 
$$
[n]_r=\frac{n!}{(n-r)!}=n(n-1)\ldots(n-r+1)
$$
The (positive) Stirling numbers of the second kind  are then defined by
$$
x^n=\sum_{k=0}^n S(n,k) [x]_k
$$

We also use the notations $X:=x(1+x)$, $y:=\frac{x}{1+x}$ and $m_r(x):=x^r$ for monomials.\\
$\PP$ (resp. $\PP_n$ ) denotes the space of polynomials (resp. of polynomials of degree at most $n$).
%-----------------------------------------------------------------------------------------------------------------------------------
%------------------------------------2. BASKAKOV OPERATOR AS DIFF OPERATOR------------------------------
%-----------------------------------------------------------------------------------------------------------------------------------

\section{Baskakov operator and its inverse as differential operators}
Let us denote 
$$
\VV_n:=\{v_{k,n}(x):= {n+k-1 \choose k} x^k (1+x)^{-(n+k)},\, k\ge 0\}
$$
the set of basic functions of the Baskakov operators
$$
\mathcal{V}_nf(x):=\sum_{k\ge 0} f_k \, v_{k,n}(x),\quad f_k:=f\left(\frac kn\right)
$$
The images of (Newton) polynomials are monomials of the same degree  (\cite{MM2}, Lemma 1.1)
$$
\nu_{r,n}(x):=\prod_{i=0}^{r-1}\left(x-\frac in\right) \;\;  \Longrightarrow\;\;  \mathcal{V}_n\nu_r(x)=\lambda_{n,r} m_r(x)
$$
where
$$
\lambda_{n,r}=\prod_{i=0}^{r-1} \left(1+\frac in\right)=\frac{(n)_r}{n^r} ,\qquad m_r(x):=x^r.
$$
Therefore $\mathcal{V}_n$ is degree-preserving for all $r\ge 0$. One can write
$$
\mathcal{V}_n\nu_r(x)=\sum_{k\ge r}\nu_r \left(\frac kn\right) v_{k,n}(x)=n^{-r}(1+x)^{-n}y^r\sum_{k\ge r} \frac{(n+k-1)!}{(k-r)!\,(n-1)!}y^{k-r}
$$
With the change of index $k=j+r$, the series can be written
$$
\sum_{j\ge 0} \frac{(n+j+r-1)!}{j!\, (n-1)!}y^j=\sum_{j\ge 0} u_j
$$
As $\lim_{j\to +\infty}\frac{u_{j+1}}{u_j}=\lim_{j\to +\infty} \frac{j+n+r}{j+1} y=y<1$, 
the series is always convergent. Moreover, $\mathcal{V}_n$ is exact on $\PP_1$ since 
$\mathcal{V}_n\nu_0(x)=m_0(x)=1$ and $\mathcal{V}_n\nu_1(x)=m_1(x)=x$.

%----------------------------------------2.1.V_n and its inverse on \PP-------------------------------------------

\subsection{$\mathcal{V}_n$ and its inverse on $\PP$}

Therefore $\mathcal{V}_n$ is an isomorphism of the space $\PP$ of polynomials and both $\mathcal{V}_n$ and $\mathcal{U}_n:=\mathcal{V}_n^{-1}$ can be written as differential operators on this space (here $D=d/dx$):
$$
\mathcal{V}_n=\sum_{r\ge 0}\theta_r^{(n)} D^r \qquad \mathcal{U}_n=\sum_{r\ge 0} \eta_r^{(n)}D^r.
$$
In the following, we often omit the upper index $n$ when the latter is fixed.

%-------------------------------------2.1.1.Computation of  polynomials theta------------------------------

\subsubsection{Computation of  polynomials $\theta_p(x)$}
We begin with computing $w_{n,p}:=\mathcal{V}_nm_p$, then solving the 
 system of linear equations \\ 
 $\sum_{r=2}^p\frac{p!}{(p-r)!}\theta_r=w_{n,p}-m_p$ in the unknowns $\theta_r$, with right-hand side the  vector with components $\mathcal{V}_nm_p-m_p$.
As we already know  that
$$
m_r(x)=\sum_{k\ge 0}\frac{[k]_r}{(n)_r} v_{k,n}(x) \quad {\rm and}\quad k^p=\sum_{r=0}^p S(p,r) [k]_r
$$
 where the $S(p,r)$ are the positive Stirling numbers of the second kind, we deduce immediately 
$$
w_{p,n}(x):=\mathcal{V}_nm_p(x)=\frac{1}{n^p}\sum_{k\ge 0} k^p v_{k,n}(x)
$$
$$
=\frac{1}{n^p}\sum_{r=0}^p S(p,r) \sum_{k\ge 0}[k]_r v_{k,n}(x)=\frac{1}{n^p}\sum_{r=0}^p S(p,r)(n)_rm_r(x)
$$
For example
$$
w_0(x)=1, \quad w_1(x)=x,  \quad w_2(x)=\frac{1}{n^2}\sum_{r=0}^2 S(2,r)(n)_rx^r=\frac{1}{n}(x+(n+1)x^2)
$$
Then $w_2(x)=m_2(x)+2\theta_2(x)$ and $\theta_2(x)=\frac12(w_2(x)-x^2)=\frac{1}{2n}x(1+x)$. Similarly
$$
w_{3}(x):=\frac{1}{n^3}\sum_{r=0}^3 S(3,r)(n)_r m_r(x)=\frac{1}{n^3}(nx+3n(n+1)x^2+n(n+1)(n+2)x^3)
$$
thus
$$
w_{3}(x)-m_3(x)=\frac{1}{n^2}(x+3(n+1)x^2+(3n+2)x^3)
$$
and $6x\theta_2+6\theta_3=w_3-m_3$ gives $\theta_3(x)=\frac{1}{6n^2}x(x+1)(2x+1)$.\\

In the general case, we get a linear system in $\theta_2,\theta_3,\ldots,\theta_r$ whose general equation is the following
$$
\frac{x^{r-2}}{(r-2)!}\theta_2+\frac{x^{r-3}}{(r-3)!}\theta_3+\ldots+ \theta_r=\pi_r(x):=\frac{1}{r!}(w_r(x)-m_r(x))
$$
By inverting the matrix, one obtains
$$
\theta_r(x)=\pi_r(x)-x\pi_{r-1}(x)+\frac{x^2}{2!}\pi_{r-2}(x)+\ldots+(-1)^r\frac{x^{r-2}}{(r-2)!}\pi_2(x)=\sum_{k=0}^{r-2} (-1)^k\frac{x^{k}}{k!}\pi_{r-k}(x)
$$
%--------------------------------------Example 1 : computation of theta_4-----------------------
{\bf Example 1: computation of $\theta_4$}\\
From the expansion $\theta_4(x)=\pi_4(x)-x\pi_{3}(x)+\frac{x^2}{2!}\pi_{2}(x)$
where $\pi_2(x)=\frac{1}{2n}x(1+x)$, \\
$\pi_3(x)=\frac{1}{6n^2}(x+3(n+1)x^2+(3n+2)x^3)$ and
$\pi_4(x)=\frac{1}{4!}(w_4(x)-m_4(x))$, with
$$
w_4(x)=\frac{1}{n^4}\sum_{r=0}^4 S(4,r)(n)_r m_r(x)=\frac{1}{n^4}(nx+7n(n+1)x^2+6(n)_3x^3+(n)_4x^4)
$$
we deduce
$$
\pi_4(x)=\frac{1}{4!}(w_4(x)-m_4(x))=\frac{1}{4! n^3}(x+7(n+1)x^2+6(n+1)(n+2)x^3+(6n^2+11n+6)x^4)
$$
and finally we obtain the right polynomial
$$
\theta_4(x)=\frac{X}{24n^3}(1+3(n+2)X)
$$
%--------------------------------------Example 2 : computation of theta_5-----------------------
{\bf Example 2: computation of $\theta_5$}\\
Similarly, from the expansion $\theta_5(x)=\pi_5(x)-x\pi_{4}(x)+\frac{x^2}{2!}\pi_{3}(x)-\frac{x^3}{3!}\pi_{2}(x)$ where
$$
w_5(x)=\frac{1}{n^5}\sum_{r=0}^5 S(5,r)(n)_r m_r(x)=\frac{1}{n^5}(nx+15(n)_2x^2+25(n)_3x^3+10(n)_4x^4+(n)_5x^5)
$$
from which we derive the expression of  $\pi_5(x)=\frac{1}{5!}(w_5(x)-m_5(x))$:
$$
\pi_5(x)=\frac{1}{5! n^4}(x+15(n+1)x^2+25(n+1)_2x^3+10(n+1)_3x^4+(10n^3+35n^2+50n+24)x^5),
$$
 we get
$$
\theta_5(x)=\frac{x}{5! n^4}(1+(10n+15)x+(40n+50)x^2+(50n+60)x^3+(20n+24)x^4)
$$
which can be factorized as
$$
\theta_5(x)=\frac{X}{5! n^4}(2x+1)(1+(10n+12)X)
$$
%---------------------------------2.1.2.Computation of polynomials eta---------------------------------------

\subsection{Computation of  polynomials  $\eta_r^{(n)}$}

In \cite{MM2}, one finds
$$
\eta_2(x)=-\frac12 \frac{X}{n+1}, \quad \eta_3(x)=\frac13 \frac{(1+2x)X}{(n+1)(n+2)},\quad 
\eta_4(x)=-\frac18 \frac{X(2-(n-6)X)}{(n+1)(n+2)(n+3)}
$$
From $\mathcal{V}_n ([nx]_r)=(n)_r m_r(x)$, we deduce $\mathcal{U}_n m_r(x)=[nx]_r/(n)_r$ and, using the expansion $\mathcal{U}_n=\sum_{r\ge 0} \eta_r^{(n)}D^r$, 
we see that  $\eta=(\eta_k)$ is solution of the linear system 
$$
\eta_k+\sum_{i=0}^{k-1}\frac{x^{k-i}}{(k-i)!}\eta_i=\rho_k:=\frac{1}{k!} \frac{[nx]_k}{(n)_k}
$$
whose solution can be written under the form
$$
\eta_k=\sum_{j=0}^k (-1)^j \frac{x^j}{j!} \rho_{k-j}.
$$

%--------------------------------------Verification: computation of eta_4-----------------------
{\bf Example: computation of $\eta_4$}\\
We have
$$
\eta_4=\sum_{j=0}^4 (-1)^j \frac{x^j}{j!} \rho_{4-j},\qquad \rho_k:=\frac{1}{k!} \frac{[nx]_k}{(n)_k}
$$
for $k=0,\ldots,4$, which gives
$$
\rho_0=1, \qquad  \rho_1=x,\qquad \rho_2=\frac12\frac{nx(nx-1)}{n(n+1)}
$$
$$
\rho_3=\frac16\frac{nx(nx-1)(nx-2)}{n(n+1)(n+2)},\qquad \rho_4=\frac{1}{24}\frac{nx(nx-1)(nx-2)(nx-3)}{n(n+1)(n+2)(n+3)}
$$
Therefore
$$
\eta_4=\rho_{4}-x\rho_{3}+\frac{x^2}{2}\rho_{2}-\frac{x^3}{6}\rho_{1}+\frac{x^4}{24}
$$
Substituting the polynomials $\rho_k$ by their expressions and factorizing, we get 
$$
\eta_4(x)=-\frac18 \frac{X(2-(n-6)X)}{(n+1)(n+2)(n+3)}$$
as desired. However, the polynomials $\theta_r$ and $\eta_r$ are more easily calculated by using the recurrence formulas of the next theorem.

%--------------------------Theorem 1: Recurrence on polynomials \beta and \alpha----------------------

\subsection{Recurrence on polynomials $\beta$ and $\alpha$}
{\bf Theorem 1.} {\em The polynomials $\theta_r$ and $\eta_r$ satisfy the following recurrence relations:}
$$
n(r+1)\theta_{r+1}(x)=X(D\theta_r(x)+\theta_{r-1}(x)),\quad \theta_0=1,\theta_1=0.
$$
$$
(n+r)(r+1)\eta_{r+1}(x)=-r(1+2x)\eta_r(x)-X\eta_{r-1}(x),\quad \eta_0=1,\eta_1=0.
$$
%---------------------------------------------Proof recurrence theta------------------------------------
{\sl Proof.} For the polynomials $\theta_r$, the recurrence relation is already given in Lemma 1.1 of  \cite{MM2}. For the polynomials $\eta_r$, the proof is rather technical, so we omit some details in calculations.
The idea is to show that $E_r=\sum_{k=2}^r \frac{x^{r-k}}{(r-k)!} A_k=0$ where
$$
A_k:=k(n+k-1) \eta_k+(k-1)(2x+1)\eta_{k-1}+X\eta_{k-2},\quad k\ge 2.
$$
Moreover, as we have seen above, $\eta=(\eta_k)$ is solution of the linear system 
$$
\eta_k+\sum_{i=0}^{k-1}\frac{x^{k-i}}{(k-i)!}\eta_i=\rho_k:=\frac{1}{k!} \frac{[nx]_k}{(n)_k}
$$
and it can thus be written under the form
$$
\eta_k=\sum_{j=0}^k (-1)^j \frac{x^j}{j!} \rho_{k-j}
$$
Substituting the polynomials $\eta_j$ for their expansion in terms of polynomials $\rho_j$ in the expression
$$
E_r=\sum_{k=2}^r \frac{x^{r-k}}{(r-k)!} \left\{k(n+k-1)\eta_k+(k-1)(2x+1)\eta_{k-1}+X\eta_{k-2}\right\}
$$
we obtain
$$
E_r=\sum_{k=0}^{r-2} \frac{x^{k}}{k!}\left\{a_{k}\rho_{r-k}+b_{k}(2x+1)\rho_{r-k-1}+Xc_{k}\rho_{r-k-2} \right\}
$$
where
$$
a_k:=\sum_{i=0}^k (-1)^i {k\choose i}(r-k+i)(n+r-k+i-1)
$$
$$
b_k:=\sum_{i=0}^k (-1)^i {k\choose i}(r-1-k+i),\quad 
c_k:=\sum_{i=0}^k (-1)^i {k\choose i}
$$
For $k=0,1,2$, we get respectively $c_0=1,b_0=r-1,a_0=r(n+r-1)$, $c_1=0,b_1=-1$, $a_1=-(n+2r-2)$ and $c_2=b_2=0$, $a_2=2$.
For $k\ge 3$, it is easy to prove that $a_k=b_k=c_k=0$, 
therefore it only remains
$$
E_r=r(n+r-1)\rho_r+(r-1)(2x+1)\rho_{r-1}+X\rho_{r-2}-(n+2r-2)x\rho_{r-1}-x(2x+1)\rho_{r-2}+x^2\rho_{r-2}
$$
which, using the definition of polynomials $\rho_r$ and $\rho_{r-1}$ and after simplification, gives the desired result:
$$
E_r=r(n+r-1)\rho_r-(nx-(r-1))\rho_{r-1}=0
$$
Now, it is easy to check that $A_2=0$. Then, by induction on $r$, assuming that  $A_k=0$ for $k=2,\ldots r-1$ and using the fact  that $E_r=0$, we immediately deduce that $A_r=0$, q.e.d.
$\blacksquare$
%------------------------------------2.3. ASYMPTOTIC BEHAVIOUR OF THETA-ETA-----------------------------

\subsection{Asymptotic behaviour of $\theta$ and $\eta$}
From the above results, we deduce, when $n\to \infty$:
$$
\lim n\theta_2^{(n)}=\frac12X,\quad  \lim n^2 \theta_3^{(n)}=\frac16 (1+2x)X,\quad \lim n^2 \theta_4^{(n)}=\frac18X^2,\quad \lim n^3\theta_5^{(n)}=\frac{1}{12}(1+2x)X^2
$$
$$
\lim n\eta_2^{(n)}=-\frac12X,\quad  \lim n^2 \eta_3^{(n)}=\frac13(1+2x)X,\quad \lim n^2 \eta_4^{(n)}=\frac18X^2,\quad \lim n^3\eta_5^{(n)}=\frac16(1+2x)X^2
$$
More generally, using the recurrence relations of Theorem 1, one obtains, by induction on $r$:\\

%----------------------------------------------------Theorem 2------------------------------------
{\bf Theorem 2.}{\em The following limits hold, when $n\to +\infty$}
$$
\lim n^r\theta_{2r-1}^{(n)}=\bar\theta_{2r-1}=\frac{1}{3\cdot 2^{r-1}(r-2)!}(1+2x)X^{r-1},\qquad \lim n^r \theta_{2r}^{(n)}=\bar\theta_{2r}=\frac{1}{2^rr!}X^r
$$
$$
\lim n^r\eta_{2r-1}^{(n)}=\bar\eta_{2r-1}=\frac{(-1)^r}{3\cdot 2^{r-2}(r-2)!}(1+2x)X^{r-1},\qquad \lim n^r \eta_{2r}^{(n)}=\bar\eta_{2r}=\frac{(-1)^r}{2^{r} r!}X^r
$$
{\em Proof.} The result being already true for $r=1,2$, we assume the limits hold for $r$ and we prove them for $r+1$. We only give the proof for polynomials $\eta_{2r+1}^{(n)}$ and $\eta_{2r+2}^{(n)}$, the one for $\theta_{2r+1}^{(n)}$ and $\theta_{2r+2}^{(n)}$ being similar. We start respectively from
$$
(n+2r)n^r (2r+1)\eta^{(n)}_{2r+1}(x)=-2r(1+2x)n^r\eta^{(n)}_{2r}(x)-Xn^r\eta^{(n)}_{2r-1}(x)
$$
$$
(n+2r+1)n^r(2r+2)\eta^{(n)}_{2r+2}(x)=-(2r+1)(1+2x)n^r\eta^{(n)}_{2r+1}(x)-Xn^r\eta^{(n)}_{2r}(x)
$$
Taking the limits, which exist, when $n\to +\infty$, we get
$$
(2r+1)\bar\eta_{2r+1}=-2r(1+2x)\bar\eta_{2r}-X\bar\eta_{2r-1}
$$
$$
(2r+2)\bar\eta_{2r+2}=-X\bar\eta_{2r}
$$
whence the respective polynomials:
$$
\bar\eta_{2r+1}=-\frac{2r}{2r+1}(1+2x)\bar\eta_{2r}-\frac{X}{2r+1}\bar\eta_{2r-1}
$$
$$
=(-1)^{r+1}(1+2x)X^r\frac{1}{2r+1} \left(\frac{1}{2^{r-1} (r-1)!}+\frac{1}{3\cdot 2^{r-2}(r-2)!}\right)
$$
$$
=(-1)^{r+1}(1+2x)X^r\frac{1}{3\cdot 2^{r-1} (r-1)!}
$$
$$
\bar\eta_{2r+2}=-\frac{1}{(2r+2)}X\bar\eta_{2r}=\frac{(-1)^{r+1}}{2^{r+1} (r+1)!}X^{r+1}
$$
which completes the proof. $\blacksquare$

%-----------------------------------------------------------------------------------------------------------------------------------
%------------------------------------3. BASKAKOV QIs: NORMS & CONVERGENCE--------------------------------
%-----------------------------------------------------------------------------------------------------------------------------------

\section{Baskakov quasi-interpolants: norm and convergence}

Baskakov left quasi-interpolants are defined  \cite{MM2}  by
$$
\mathcal{V}_n^{(r)} f(x)=\mathcal{U}_n^{(r)}\mathcal{V}_n f(x)=\sum_{k=0}^r \eta_k^{(n)}(x)D^k \mathcal{V}_n f(x).
$$
(However in that paper the authors do not give the recurrence relation of Theorem 1 on the coefficients $\eta$).
Baskakov right quasi-interpolants are defined by
$$
\mathcal{V}_n^{[r]} f(x)=\mathcal{V}_n\mathcal{U}_n^{(r)}f(x)=\sum_{k=0}^r \mathcal{V}_n (\eta_k^{(n)}D^k f)(x).
$$
They seem to be less interesting than the left ones since they need the knowledge of derivatives of the function $f$ to be approximated. Therefore, we do not study them in the present paper. In the rest of the paper, the QIs are thus left QIs.

%-----------------------------------------Norm of Baskakov QIs------------------------------

\subsection{Norms of Baskakov quasi-interpolants}
The (left) Baskakov QIs can also be represented in the {\em quasi-Lagrange form}
$$
\mathcal{V}_n^{(r)} f(x)=\sum_{j\ge 0} f\left(\frac jn\right) v_{j,n}^{(r)}(x),
$$
where the {\em quasi-Lagrange basic functions} are defined by
$$
v_{j,n}^{(r)}(x):=\sum_{k=0}^r \eta_k^{(n)}(x)D^k  v_{j,n}(x)
$$ 
This form allows to express the infinite norm of the operator $\mathcal{V}_n^{(r)} $ as the max norm of the associated Lebesgue function
$$
\Vert \mathcal{V}_n^{(r)}  \Vert_\infty=| \Lambda^{(r)}_n |_\infty=\max_{x\ge 0} \left(\sum_{j\ge 0}\, |v_{j,n}^{(r)} (x)| \right)
$$
The following result is proved in \cite{MM2} (Lemma 2.3):\\

%----------------------------------------------------Theorem 3------------------------------------------------
{\bf Theorem 3.} {\em For all $r\ge 0$, there exists a constant $C_r>0$ such that}
$$
\Vert \mathcal{V}_n^{(r)}  \Vert_\infty \le C_r \quad \forall n\ge r
$$
The following table gives the first values of the norms and estimations of the expected values of $C_r$.
\begin{center}
\begin{tabular}{|c|c|c|c|c|c|c|c|c|c|}
\hline
$n\backslash r$&2&3&4&5&6&7&8&9\\
\hline
8&1.10&1.32&1.72&2.34&3.30&4.90&7.50&\\
\hline
16&1.12&1.34&1.77&2.44&3.50&5.20&8.00&13.0\\
\hline
24&1.128&1.345&1.79&2.46&3.54&5.31&8.26&13.3\\
\hline
32&1.131&1.350&1.80&2.48&3.58&5.38&8.38&13.5\\
\hline
40&1.133&1.352&1.80&2.50&3.60&5.42&8.46&13.7\\
\hline
48&1.134&1.353&1.80&2.50&3.62&5.44&8.51&13.8\\
\hline
$C_r$&1.15&1.4&1.9&2.6&3.7&5.6&8.8&14\\
\hline
\end{tabular}
\end{center}
At first sight, it can be observed that the norms have reasonable values for $1\le r\le 7$ and much higher ones when $r\ge 8$. 
%----------------------------------------3.1. Asymptotic convergence order-------------------------------------
\subsection{Asymptotic convergence order}
From Theorems 1 and 2, we may immediately deduce the asymptotic convergence order of QIs for polynomials.\\

%----------------------------------------------------Theorem 4------------------------------------
{\bf Theorem 4.} {\em Let $p$ be a polynomial of degree  $m\ge 2r+1$, then}
$$
\lim n^{r+1}(p(x)-\mathcal{V}_n^{(2r)}p(x))=\bar \eta_{2r+1} D^{2r+1}p(x)+\bar \eta_{2r+2} D^{2r+2}p(x)
$$
$$
\lim n^{r+1}(p(x)-\mathcal{V}_n^{(2r+1)}p(x))=\bar \eta_{2r+2} D^{2r+2}p(x)
$$
{\sl Proof.} For $p\in \PP_{m}$, with $m\ge 2r+1$, one can write
$$
p(x)-\mathcal{V}_n^{(2r)}p(x)=(\mathcal{U}_n-\mathcal{U}_n^{(r)})\mathcal{V}_np(x)=\sum_{k=2r+1}^m  \eta_k^{(n)}(x)D^k \mathcal{V}_np(x)
$$
Remind that $\lim n^{r+1} \eta_{2r+1}^{(n)}=\frac{(-1)^{r+1}}{3\cdot 2^{r-1}(r-1)!}(1+2x)X^r$, $\lim  n^{r+1} \eta_{2r+2}^{(n)}=\frac{(-1)^{r+1}}{2^{r+1} (r+1)!} X^{r+1}$ and $\lim n^{r+1}\eta_k^{(n)}=0$ for $k\ge 2r+3$. Moreover,
$\lim D^k \mathcal{V}_np(x)=D^kp(x)$, therefore
$$
\lim n^{r+1}\sum_{k=2r+1}^m  \eta_k^{(n)}(x)D^k \mathcal{V}_np(x)=\lim n^{r+1} \eta_{2r+1}^{(n)}D^{2r+1}p(x)+\lim n^{r+1} \eta_{2r+2}^{(n)}D^{2r+2}p(x).
$$
The proof is similar for $\mathcal{V}_n^{(2r+1)}$. $\blacksquare$\\

In other terms
$$
p(x)-\mathcal{V}_n^{(2r)}p(x)=O(n^{-(r+1)})\qquad p(x)-\mathcal{V}_n^{(2r+1)}p(x)=O(n^{-(r+1)})
$$
This result can be extended, via Taylor's formula, to functions which are differentiable enough.\\

%------------------------------------------Theorem 5: VORONOVSKAYA -------------------------------------------------

{\bf Theorem 5 (Voronovskaya)} {\em Let $f\in C^{2r+4}(\RR_+)$ be a function with $D^{2r+4}f$ bounded, then for all $x\ge 0$}
$$
\lim n^{r+1}(f(x)-\mathcal{V}_n^{(2r)}f(x))=\bar \eta_{2r+1} D^{2r+1}f(x)+\bar \eta_{2r+2} D^{2r+2}f(x)
$$
$$
\lim n^{r+1}(f(x)-\mathcal{V}_n^{(2r+1)}f(x))=\bar \eta_{2r+2} D^{2r+2}f(x)
$$

%-------------------------------------------------PROOF for any r-----------------------------------------------------

{\sl Proof}. We only give the proof for the second limit, the one for the first limit being similar.
$$
\lim n^{r+1}(f(x)-\mathcal{V}_n^{(2r+1}f(x))=\bar \eta_{2r+2} D^{2r+2}f(x)
$$
By Taylor's formula, for $x$ fixed and for all $t$, there exists $x_t$ between $x$ and $t$ such that 
$$
f(t)=f(x)+\sum_{k=0}^{2r+3} \frac{1}{k!} (t-x)^k D^k f(x)+\frac{1}{(2r+4)!} (t-x)^{2r+4} D^{2r+4}f(x_t)
$$
Therefore, as $\mathcal{V}_n^{(2r+1)}$ is exact on $\PP_{2r+1}$:
$$
\mathcal{V}_n^{(2r+1)}f(x)-f(x)= \frac{1}{(2r+2)!} D^{2r+2} f(x)\mathcal{V}_n^{(2r+1)} (t-x)^{(2r+2)}
$$
$$
+\frac{1}{(2r+3)!} D^{(2r+3)}f(x)\mathcal{V}_n^{(2r+1)} (t-x)^{(2r+3)}+\frac{1}{(2r+4)!}\mathcal{V}_n^{(2r+1)} [(t-x)^{(2r+4)}D^{(2r+4)}f(x_t)]
$$
First, we get $\mathcal{V}_n^{(2r+1)}[(t-x)^{(2r+2)}]=\mathcal{V}_n^{(2r+1)}m_{2r+2}-m_{2r+2}=-\eta_{2r+2}D^{(2r+2)}\mathcal{V}_n m_{2r+2}$ since by definition
$$
\mathcal{V}_n^{(2r+2)}m_{2r+2}=m_{2r+2}=\mathcal{V}_n^{(2r+1)}m_{2r+2}+\eta_{2r+2}(x)D^{(2r+2)}\mathcal{V}_n m_{2r+2}
$$  
Moreover, as $\lim D^{2r+2}\mathcal{V}_n m_{2r+2}=D^{2r+2}m_{2r+2}=(2r+2)!$ and $\lim n^{r+1}\eta_{2r+2}(x)=\bar \eta_{2r+2}(x)$, we obtain
$$
\lim n^{r+1} \left(\frac{1}{(2r+2)!} D^{(2r+2)} f(x)(\mathcal{V}_n^{(2r+1)}m_{2r+2}-m_{2r+2})\right)=-\bar \eta_{2r+2}(x)D^{2r+2} f(x)
$$
Second, we will prove that $\lim n^{r+1} \mathcal{V}_n^{(2r+1)} [(t-x)^{2r+3}]=0$. We first obtain
$$
\mathcal{V}_n^{(2r+1)}[(t-x)^{2r+3}]=(\mathcal{V}_n^{(2r+1)}m_{2r+3}-m_{2r+3})-(2r+3)x(\mathcal{V}_n^{(2r+1)}m_{2r+2}-m_{2r+2})
$$
We already know that $\lim n^{r+1} (\mathcal{V}_n^{(2r+1)}m_{2r+2}-m_{2r+2})=-(2r+2)! \bar \eta_{2r+2}(x)$. Then, using the fact that $ \mathcal{V}_n^{(2r+3)} m_{2r+3}=m_{2r+3}$, we get
$$
n^{r+1}(\mathcal{V}_n^{(2r+1)}m_{2r+3}-m_{2r+3})=-n^{r+1}\eta_{2r+2}(x) D^{2r+2}\mathcal{V}_n m_{2r+3}-n^{r+1}\eta_{2r+3}(x) D^{2r+3}\mathcal{V}_n m_{2r+3}
$$
As $\lim D^{2r+2}\mathcal{V}_n m_{2r+3}=D^{2r+2}m_{2r+3}=(2r+3)! x$, $\lim D^{2r+3}\mathcal{V}_n m_{2r+3}=D^{2r+3}m_{2r+3}=(2r+3)! $, $\lim n^{r+1}\eta_{2r+2}(x)=\bar \eta_{2r+2}(x)$ and $\lim n^{r+1}\eta_{2r+3}(x)=0$, we obtain
$$
\lim n^{r+1}(\mathcal{V}_n^{(2r+1)}m_{2r+3}-m_{2r+3})=-(2r+3)! \,x\, \bar\eta_{2r+2}(x)
$$
and finally
$$
\lim n^{r+1} \mathcal{V}_n^{(2r+1)} [(t-x)^{2r+3}]=-(2r+3)! \,x\, \bar\eta_{2r+2}(x)+(2r+3)\,x\, (2r+2)!\, \bar \eta_{2r+2}(x)=0
$$

It remains to prove that $\lim n^{r+1} R(x)=0$ , where $R$ is the quantity
$$
R:=\mathcal{V}_n^{(2r+1)} [(t-x)^{2r+4}D^{2r+4}f(x_t)]=\sum_{j\ge 0} f^{(2r+4)}(x_j)\left( \frac jn-x\right)^{2r+4} v_{n,j}^{(2r+1)}(x)
$$
where the quasi-Lagrange functions are defined by
$$
v_{n,j}^{(2r+1)}(x):=v_{n,j}(x)+\sum_{k=2}^{2r+1}\eta_k(x)D^kv_{n,j}(x)
$$
We have seen above  that
$$
X^2D^2v_{n,j}(x)=p_2(x)v_{n,j}(x),\quad {\rm where} \;\; |p_2(x)|\approx n^2x^2
$$
$$
X^3 D^3v_{n,j}(x)=p_3(x)v_{n,j}(x), \quad {\rm where}\;\;  |p_3(x)|\approx n^3 x^3
$$
More generally, it is easy to prove that 
$$
X^k D^kv_{n,j}(x)=p_k(x)v_{n,j}(x), {\rm where} \;\; |p_k(x)|\sim n^{k} x^{k}
$$
We then obtain
$$
|v_{n,j}^{(2r+1)}(x)|\le v_{n,j}(x)\left(1+\sum_{k=2}^{2r+1}X^{-k}|\eta_k(x)||p_k(x)|\right)
$$
From  Theorem 2, we deduce
$$
|\eta_{k}(x)|\sim n^{-l} \bar \eta_k(x),\quad {\rm for}\;\; k=2l-1\;\;{\rm or}\;\; k=2l
$$
thus
$$
|\eta_{2l-1}(x)||p_{2l-1}(x)|\sim  n^{l-1}x^{2l-1}\bar \eta_{2l-1}(x),\quad
|\eta_{2l}(x)||p_{2l}(x)|\sim  n^{l}x^{2l}\bar \eta_{2l}(x),\
$$
Therefore the parenthesis in the upper bound on $|v_{n,j}^{(2r+1)}(x)|$ is equivalent, when $n\to\infty$, to
$$
X^{-2r}|\eta_{2r}(x)||p_{2r}(x)|+X^{-(2r+1)}|\eta_{2r+1}(x)||p_{2r+1}(x)|\sim n^{r}\left(X^{-2r}\bar \eta_{2r}(x)+X^{-(2r+1)}\bar \eta_{2r+1}(x)\right)
$$
Let us introduce the two following subsets of indices, for a given $a>0$:
$$
J_1:=\{ j\in \NN\, :\, |j/n-x|<n^{-a}\},\quad  J_2:=\{ j\in \NN\, :\, |j/n-x|\ge n^{-a}\}
$$
We can write
$$
R\le \Vert f^{(2r+4)}\Vert_\infty\left( \sum_{j\in J_1} \left( \frac jn-x\right)^{2r+4} |v_{n,j}^{(2r+1)}(x)|+\sum_{j\in J_2}\left( \frac jn-x\right)^{2r+4} |v_{n,j}^{(2r+1)}(x)|Ê\right)
$$
In $R_1:=\sum_{j\in J_1} \left( \frac jn-x\right)^{2r+4} |v_{n,j}^{(2r+1)}(x)|$, we have $ \left| \frac jn-x\right|^{2r+4}<n^{-(2r+4)a}$ \\
and $\sum_{j\in \ZZ} |v_{n,j}^{(2r+1)}(x)|\le \Vert \mathcal{V}_n^{(2r+1)}\Vert_\infty\le C_{2r+1}$,  whence $n^{r+1}R_1\le C_{2r+1}n^{r+1-2(r+2)a}$ and $\lim n^{r+1} R_1=0$ if we take $a>\frac{r+1}{2(r+2)}$.

The  second sum is equivalent, up to the rational factor  $X^{-2r}\bar \eta_{2r}(x)+X^{-(2r+1)}\bar \eta_{2r+1}(x)$, to 
$$
R_2=n^{r}\sum_{j\in J_2} \left( \frac jn-x\right)^{2r+4} v_{n,j}(x)
$$
and we have to prove that $\lim n^{r+1} R_2=\lim n^{2r+1}\sum_{j\in J_2} \left( \frac jn-x\right)^{2r+4} v_{n,j}(x)=0$.

As $j\in J_2$, $|j/n-x|\ge n^{-a}$ implies $|j-nx|^{2p}\ge n^{2p(1-a)}$, thus $1\le n^{2p(a-1)} |j-nx|^{2p}$, whence
$$
\sum_{j\in J_2} (j-nx)^{2r+4} v_{n,j}(x)\le n^{2p(a-1)} \sum_{j\in J_2} (j-nx)^{2r+2p+4} v_{n,j}(x)
$$
From the recurrence relation of Theorem 1, we already know that $\theta_{2r+2p+4}(x)\le C n^{-(r+p+2)}$, thus
$$
\sum_{j\in \ZZ} (j-nx)^{2r+2p+4} v_{n,j}(x)=n^{2r+2r+2p+4} \sum_{j\in \ZZ} (j/n-x)^{2r+2p+4} v_{n,j}(x)\le Cn^{r+p+2}
$$
This gives
$$
 n^{2r+1}\sum_{j\in J_2} \left( \frac jn-x\right)^{2r+4} v_{n,j}(x)\le Cn^{-3}n^{2p(a-1)}n^{r+p+2}=Cn^{2pa+r-p-1}
$$
and the limit will be equal to 0 provided $a<\frac{p+1-r}{2p}$. Finally, we must have the following bounds on $a$:
$$
\frac{r+1}{2(r+2)}<a<\frac{p+1-r}{2p}
$$
This is possible if $\frac{r+1}{r+2}<\frac{p+1-r}{p}$, i.e. if we choose $p>(r-1)(r+2)$.  $\blacksquare$

%-----------------------------------------------------------------------------------------------------------------------------------
%----------------------------------4. COMPUTATION OF BASKAKOV QIs----------------------------------------------
%-----------------------------------------------------------------------------------------------------------------------------------
\section{Computation of Baskakov quasi-interpolants}
We give more explicit expressions of  Baskakov quasi-interpolants $\mathcal{V}_n^{(r)}f(x)$  for the first orders $1\le r\le 11$, which can be useful for their practical evaluation. The general case can be treated in the same way.
%------------------------------------------4.1. Baskakov operator---------------------------
\subsection{Baskakov operator}
The problem is to compute, for a given $N$ large enough,
$$
\mathcal{V}_nf(x):=\sum_{k=0}^N f_k v_{k,n}(x),\quad v_{k,n}(x):= {n+k-1 \choose k} x^k (1+x)^{-(n+k)}
$$
Recalling the notations $f_k:=f\left(\frac kn\right)$ and  $y=\frac{x}{1+x}$, one gets
$$
\mathcal{V}_nf(x):=(1+x)^{-n}\sum_{k=0}^N f_k {n+k-1 \choose k} y^k=(1+x)^{-n}P_N(y)
$$ 
As $P_N(y)$ is a polynomial of degree at most $N$ in the variable $y$, it can be evaluated in $O(N)$ operations.\\

For the results below, we need the following  lemma whose proof is technical, but straightforward\\

{\bf Lemma.} {\em (i) For given pairs $(p,j)$ and $(n,k)$, with $0\le j\le p$ and $0\le k\le n$, there holds}
$$
x^j(1+x)^{p-j}v_{k,n}(x)=\omega_{p,j}(n,k) v_{k+j,n-p}(x)
$$
{\em where}
$$
\omega_{p,j}(n,k):=\frac{(n+k-1)\ldots(n+k-p+j)(k+1)\ldots(k+j)}{(n-1)\ldots(n-p)}
$$
{\em (ii)  For all $p\ge 1$, denoting as usual $\Delta u_k:=u_{k+1}-u_k$, there holds}
$$
D^p v_{k,n}=(-1)^p (n)_p\, \Delta^p v_{k-p,n+p}
$$
{\em therefore}
$$
D^p \mathcal{V}_nf(x)=(n)_p\sum_{k\ge 0} (\Delta^p f_k) \, v_{k,n+p}(x)
$$
%------------------------------------------7.2. Baskakov QI : r=2---------------------------
\subsection{Baskakov QI: $r=2$}
We have to compute
$$
\mathcal{V}_n^{(2)}f(x):=\mathcal{V}_nf(x)+\eta_2^{(n)}(x)D^2\mathcal{V}_nf(x),
$$
with $\eta_2^{(n)}(x)=-\frac{X}{2(n+1)}$ and
$$
D^2\mathcal{V}_nf(x)=n(n+1)\sum_{k=0}^{N-2} (\Delta^2 f_k) \,v_{k,n+2}(x)
$$
Using the above lemma
$$
x(1+x) v_{k,n+2}(x)= {n+k+1 \choose k} x^{k+1} (1+x)^{-(n+k+1)},
$$
we obtain
$$
\eta_2^{(n)}(x)D^2\mathcal{V}_nf(x)=-\frac{n}{2} (1+x)^{-n}y \sum_{k\ge 0} {n+k+1 \choose k} (\Delta^2 f_{k}) y^{k}
$$

%------------------------------------------7.3. Baskakov QI : r=3----------------------------------------
\subsection{Baskakov QI: $r=3$}
We have to compute
$$
\mathcal{V}_n^{(3)}f(x):=\mathcal{V}_n^{(2)}f(x)+\frac{(1+2x)X}{3(n+1)(n+2)}
D^3\mathcal{V}_nf(x)
$$
As $(1+2x)X=x(1+x)^2+x^2(1+x)$ and
$$
D^3\mathcal{V}_nf(x)=n(n+1)(n+2)\sum_{k=0}^{N-3} (\Delta^3 f_k) \,v_{k, n+3}(x)
$$
we obtain, using $(1+2x)X=(1+x)^3y(1+y)$:
$$
\mathcal{V}_n^{(3)}f(x):=\mathcal{V}_n^{(2)}f(x)+\frac n3(1+x)^{-n}y(1+y)\sum_{k=0}^{N-3} {n+k+2 \choose k}(\Delta^3 f_k)\, y^k
$$

%------------------------------------------5.4. Baskakov QI : r=4---------------------------
\subsection{Baskakov QI: $r=4$}
In a similar way, we compute
$$
\mathcal{V}_n^{(4)}f(x):=\mathcal{V}_n^{(3)}f(x)+\eta_4^{(n)}(x)D^4\mathcal{V}_nf(x),\quad \eta_4^{(n)}(x)=-\frac{X(2-(n-6)X)}{8(n+1)(n+2)(n+3)}
$$
Using the auxiliary variable $ z=y+1/y$, the numerator of $\eta_4(x)$ can be written
$$
2X-(n-6)X^2=2x(1+x)^2-(n-2)x^2(1+x)^2+2x^3(1+x)
$$
$$
=(1+x)^{4}y^2(2/y-(n-2)+2y)=(1+x)^{4}y^2(2z-(n-2)).
$$
From 
$$
D^4\mathcal{V}_nf(x)=(n)_4\sum_{k=0}^{N-4} (\Delta^4 f_k)\,v_{k, n+4}(x)
$$
we deduce the  following expression
$$
\eta_4^{(n)}(x)D^4\mathcal{V}_nf(x)=-\frac n8(1+x)^{-n}y^2(2z-(n-2))\sum_{k\ge0}{n+k+3 \choose k} (\Delta^4 f_k)y^k
$$

%------------------------------------------5.5. Baskakov QI : r³5---------------------------
\subsection{Baskakov QI: $r\ge 5$}
In the same way, we get successively
$$
r=5: \quad \eta_5^{(n)}(x)D^5\mathcal{V}_nf(x)=\frac{4n}{5!} (1+x)^{-n}\tau_5(y) \sum_{k\ge 0} {n+k+4 \choose k}(\Delta^5 f_k)\,y^{k}
$$
$$
r=6: \quad \eta_6^{(n)}(x)D^6\mathcal{V}_nf(x)=-\frac{5n}{6!} (1+x)^{-n}\tau_6(y) \sum_{k\ge0}  {n+k+5 \choose k}(\Delta^6 f_k)\,y^{k}\pi_4(z)
$$
where
$$
\tau_5(y)=y^2(1+y)(6z-5n), \quad
\tau_6(y)=y^3(24z^2-(26n-24)z+(3n^2-34n-24))
$$
$$
r=7:\quad \eta_7^{(n)}(x)D^7\mathcal{V}_nf(x)=\frac{6n}{7!} (1+x)^{-n}\tau_7(y) \sum_{k\ge 0} {n+k+6 \choose k}(\Delta^7 f_k)\,y^{k}
$$
where
$$
\tau_7(y)=y^3(1+y)(120z^2-154nz+5(7n^2-14n-120))
$$
$$
r=8: \quad \eta_8^{(n)}(x)D^8\mathcal{V}_nf(x)=-\frac{7n}{8!} (1+x)^{-n}\tau_8(y) \sum_{k\ge 0} {n+k+7 \choose k}(\Delta^8 f_k)\,y^{k}
$$
where
$$
\tau_8(y)=y^4(720z^3-36(29n-20)z^2+4(85n^2-152n-110)z-5(3n^3-100n^2+340n+544))
$$
$$
r=9: \quad \eta_9^{(n)}(x)D^9\mathcal{V}_nf(x)=\frac{8n}{9!} (1+x)^{-n}\tau_9(y) \sum_{k\ge 0} {n+k+8 \choose k}(\Delta^9 f_k)\,y^{k}
$$
where
$$
\tau_9(y)=y^4(1+y)(7! z^3-\gamma_2z^2+\gamma_1z-\gamma_0)
$$
with
$$
\gamma_2=36(223n+120),\quad \gamma_1=4(826n^2-1197n-360),
$$
$$
\gamma_0=5(63n^3-490n^2-1152n+1152)
$$
%-----------------------------------------------------------------------------------------------------------------------------------
%-----------------4.APPROXIMATION BY BASKAKOV QIs: EXAMPLES--------------------------------------
%----------------------------------------------------------------------------------------------------------------------------------

\section{Approximation by Baskakov quasi-interpolants: examples}

In practice, one can only compute an approximation of $\mathcal{V}_n^{(r)} f(x)$ involving  $N$ terms of the series, for $N$ large enough
$$
\mathcal{V}_{n,N}^{(r)} f=\sum_{s=0}^r \eta_s^{(n)}D^s \mathcal{V}_{n,N}f\qquad \mathcal{V}_{n,N}f=\sum_{k=0}^N f_k \, v_{k,n}
$$
We often choose $N=mn$ where $m=5,6$ for example: this gives a reasonable approximation on the intervals $[0,2]$ or $[0,3]$. Of course, the results depend very much both on the function to be approximated and on the interval in which one needs this approximation.\\

{\bf Notations:} In the following tables,  the first line contains the orders of quasi-interpolants:  in column ($2r+1$), we write the error $\Vert f-\mathcal{V}_n^{(2r+1)} f \Vert_\infty$.\\

%---------------------------Example 1: f(x)=exp(-x)--------------------------
{\bf Example 1.} $f(x)=\exp(-x)$, $N=5n$, errors on the interval $[0,2]$.\\

\begin{center}
\begin{tabular}{|c|c|c|c|c|c|c|c|}
\hline
$n$&1&3&5&7&9&11\\
\hline
10&4.0(-2)&2.8(-3)&3.2(-4)&2.5(-5)&1.9(-5)&1.2(-5)\\
\hline
20&2.1(-2)&1.0(-3)&4.8(-6)&4.0(-6)&8.0(-7)&5.6(-7)\\
\hline
30&1.4(-2)&5.2(-4)&4.2(-6)&1.8(-7)&7.6(-8)&6.8(-8)\\
\hline
40&1.0(-2)&3.2(-4)&3.3(-6)&5.4(-8)&2.3(-8)&3.5(-9)\\
\hline
50&8.4(-3)&2.1(-4)&2.2(-6)&1.5(-8)&2.1(-9)&1.8(-9)\\
\hline
\end{tabular}
\end{center}
%---------------------------Example 2: f(x)=1/(1+x^2)-----------------------------
{\bf Example 2.} $f(x)=\frac{1}{1+x^2}$, $N=5n$, errors on the interval $[0,2]$.

\begin{center}
\begin{tabular}{|c|c|c|c|c|c|c|c|c|}
\hline
$n$&1&3&5&7&9&11\\
\hline
10&7.2(-2)&7.0(-3)&2.0(-3)&1.7(-3)&2.8(-4)&\\
\hline
20&3.0(-2)&2.6(-3)&2.8(-4)&1.8(-4)&5.2(-5)&1.2(-5)\\
\hline
30&2.0(-2)&1.0(-3)&1.4(-4)&3.2(-5)&1.4(-5)&2.3(-6)\\
\hline
40&1.5(-2)&8.5(-4)&8.0(-5)&7.4(-6)&4.0(-6)&8.0(-7)\\
\hline
50&1.2(-2)&5.6(-4)&4.8(-5)&2.3(-6)&1.2(-6)&2.9(-7)\\
\hline
\end{tabular}
\end{center}

%---------------------------Example 3: f(x)=exp(-x^2)-------------------------------------
{\bf Example 3.} $f(x)=\exp(-x^2)$, $N=5n$, errors on the interval $[0,2]$.

\begin{center}
\begin{tabular}{|c|c|c|c|c|c|c|}
\hline
$n$&1&3&5&7&9&11\\
\hline
10&1.0(-1)&1.7(-2)&6.4(-3)&5.0(-3)&1.3(-3)&\\
\hline
20&6.0(-2)&6.8(-3)&1.2(-3)&8.0(-4)&2.4(-4)&5.4(-5)\\
\hline
30&4.2(-2)&3.6(-3)&6.8(-4)&1.8(-4)&6.8(-5)&6.4(-6)\\
\hline
40&3.2(-2)&2.2(-3)&3.8(-4)&5.8(-5)&2.1(-5)&2.3(-6)\\
\hline
50&2.6(-2)&1.5(-3)&2.4(-4)&2.2(-5)&7.6(-6)&9.2(-7)\\
\hline
\end{tabular}
\end{center}

%---------------------------Example 4: f(x)=ln(1+x)---------------------------
{\bf Example 4.} $f(x)=\ln(1+x)$, $N=5n$, errors on the interval $[0,2]$.\\

\begin{center}
\begin{tabular}{|c|c|c|c|c|c|c|c|}
\hline
$n$&1&3&5&7&9&11\\
\hline
10&3.5(-2)&1.0(-3)&8.6(-3)&1.6(-3)&7.4(-3)&4.0(-3)\\
\hline
20&1.7(-2)&6.8(-4)&7.2(-4)&4.6(-4)&2.8(-4)&2.1(-4)\\
\hline
30&1.1(-2)&3.0(-4)&4.2(-5)&9.4(-5)&1.3(-5)&1.9(-5)\\
\hline
40&8.4(-3)&1.8(-4)&4.3(-6)&8.0(-6)&8.8(-6)&1.3(-6)\\
\hline
50&6.8(-3)&1.2(-4)&3.0(-6)&3.5(-7)&9.2(-7)&5.8(-7)\\
\hline
\end{tabular}
\end{center}

%---------------------------Example 5: f(x)=sin(6x)/(1+x^2)----------------------
{\bf Example 5.} $f(x)=\sin(6x)/(1+x^2)$, $N=5n$, errors on the interval $[0,2]$.\\
This function is oscillating and more difficult to approximate by Baskakov QIs. 

\begin{center}
\begin{tabular}{|c|c|c|c|c|c|c|}
\hline
$n$&1&3&5&7&9&11\\
\hline
10&0.64&0.44&0.32&0.30&0.24&0.23\\
\hline
20&0.43&0.28&0.18&0.15&9.6(-2)&5.0(-2)\\
\hline
30&0.36&0.20&0.11&6.8(-2)&4.0(-2)&1.2(-2)\\
\hline
40&0.30&0.16&6.8(-2)&3.2(-2)&1.8(-2)&4.7(-3)\\
\hline
50&0.23&0.13&4.6(-2)&1.5(-2)&8.0(-3)&2.2(-3)\\
\hline
\end{tabular}
\end{center}

{\bf Remarks.} \\
1) Some comments on the above computations:
\begin{itemize}
\item We are aware of the fact that this numerical study is neither complete nor sufficient  and that we should compare with other approximation methods. However, this was not the main aim of this paper.
\item For smooth and decreasing functions, the results are as expected with an improvement in the approximation order. The numerical convergence orders tend to the actual ones, predicted in Theorem 5, but we observe that  the convergence is slow for higher orders.
\item It is a little bit surprising that example 3 ($f(x)=\exp(-x^2)$) does not give as good results as example 1 ($f(x)=\exp(-x)$). Maybe this is due to the change of convexity of the former. Example 5 shows that the approximation deteriorates for an oscillating function and the observed convergence orders are not the predicted ones.
\end{itemize}
2) It is of course  possible to approximate functions on wider intervals. In that case, one has to increase the value of $N$ in order to have a good approximation. The convenient  choice of $N$ is still an open problem. \\
3) In forthcoming papers, we plan to study some techniques allowing a better approximation of series (by using e.g. convergence acceleration methods) involved in the computation of Baskakov quasi-interpolants and to give some applications to numerical analysis.

%-----------------------------------------------------------------------------------------------------------------------------------
%*******************************6. FIRST POLYNOMIALS THETA & ETA*****************************
%----------------------------------------------------------------------------------------------------------------------------------

\section{A table of the first polynomials  $\theta_r^{(n)}$ and $\eta_r^{(n)}$}

For the reader's convenience, we give tables of the first polynomials $\theta_r^{(n)}$ and $\eta_r^{(n)}$ involved in the differential representations of $\mathcal{V}_n$ and its inverse as differential operators on polynomials.
We recall that $X:=x(x+1)$ and $(n+1)_r:=(n+1)\ldots(n+k)$.\\
Using the recurrence relation of Theorem 1, Section 2, we obtain successively

%------------------------------------------6.1. Polynomials theta--------------------------------------
\subsection{Polynomials theta}
$$
\theta_6(x)=\frac{X}{6!n^5}(1+5((n+6)X+(3n^2+26n+24)X^2)
$$
$$
\theta_7(x)=\frac{X}{7!n^6}(2x+1)(1+4(14n+15)X+3(35n^2+154n+120)X^2)
$$
\smallskip
$$
\theta_8(x)=\frac{X}{8!n^7}(1+a_1X+a_2X^2+a_3X^3)
$$
$$
a_1=119n+126,\quad a_2=490n^2+2156n+1680
$$
$$
a_3=105n^3+2380n^2+7308n+5040
$$
\smallskip
$$
\theta_9(x)=\frac{X}{9!n^8}(2x+1)(1+c_1X+c_2X^2+c_3X^3)
$$
$$
c_1=246n+252,\quad c_2=1918n^2+6948n+5040
$$
$$
c_3=1260n^313216n^2+32112n+20160
$$
\smallskip
$$
\theta_{10}(x)=\frac{X}{10!n^9}(1+d_1X+d_2X^2+d_3X^3+d_4X^4)
$$ 
$$
d_1=501n+510,\quad d_2=6825n^2+24438n+17640
$$
$$
d_3=9450n^3+99120n^2+240840n+151200
$$
$$
d_4=945n^4+44100n^3+303660n^2+623376n+362880
$$
\smallskip
$$
\theta_{11}(x)=\frac{X}{11!n^{10}}(2x+1)(1+e_1X+e_2X^2+e_3X^3+e_4X^4)
$$
$$
e_1=1012n+1020,\quad e_2=22935n^2+75834n+52920
$$
$$
e_3=56980n^3+465960n^2+1013760n+604800
$$
$$
e_4=17325n^4+352660n^3+1839420n^2+3318480n+1814400
$$

%------------------------------------------6.2. Polynomials eta--------------------------------------
\subsection{Polynomials eta}
$$
\eta_5(x) := \frac{X}{30(n+1)_4}(2x+1)(6-(5n-12)X)
$$
$$
\eta_6(x) :=-\frac{x}{144(n+1)_5}(24-2(13n-60)X+(3n^2-86n+120)X^2)
$$
$$
\eta_7(x)=\frac{X}{840(n+1)_6}(2x+1)(120-(154n-480)X+(35n^2-378n+360)X^2)
$$
\smallskip
$$
\eta_8(x)=-\frac{X}{5760 (n+1)_7}(720+a_1X+a_2X^2+a_3X^3)
$$
$$
a_1=-1044n+5040, \quad a_2=340n^2-5784n+10080
$$
$$
a_3=-15n^3+1180n^2-7092n+5040.
$$
\smallskip
$$
\eta_9(x)=\frac{X}{45360(n+1)_8}(2x+1)(5040+c_1X+c_2X^2+c_3X^3)
$$
$$
c_1=-8028n+30240,\quad c_2=3304n^2-36900n+50400,
$$
$$
c_3=-315n^3+9058n^2-35928n+20160.
$$
\smallskip
$$
\eta_{10}(x)=-\frac{X}{403200 (n+1)_9}(40320+d_1X+d_2X^2+d_3X^3+d_4X^4)
$$
$$
d_1=-69264n+362880, \quad d_2=33740n^2-528912n+1088640
$$
$$
d_3=-4900n^3+199640n^2-1214880n+1209600
$$
$$
d_4=105n^4-17500n^3+273420n^2-787824n+362880
$$
\smallskip
$$
\eta_{11}(x)=-\frac{X}{3991680(n+1)_{10}}(362880+e_1X+e_2X^2+e_3X^3+e_4X^4)
$$
$$
e_1=-663696n+2903040, \quad e_2=367884n^2-4424112n+7620480
$$
$$
e_3=-70532n^3+1854072n^2-8680320n+7257600,
$$
$$
e_4=3465n^4-207284n^3+2096028n^2-4664880n+1814400.
$$
\bigskip

{\bf Acknowledgements:} I thank very much Francisco Javier Mu\~noz Delgado and the research group on Approximation of the University of Ja\'en for their kind hospitality during my stay in June 2013 during which I had the possibility of developing the present research.
%-----------------------------------------------------------------------------------------------------------------------------------
%************************************BIBLIOGRAPHY*************************************************
%-----------------------------------------------------------------------------------------------------------------------------------

\bigskip
%------------------------------------------Address-----------------------------------------
{\bf Address.} Paul Sablonni\`ere, INSA de Rennes,\\
20, Avenue des Buttes de Co\" esmes, CS 70839,\\
F-35708 Rennes Cedex 7, France.\\

{\bf E-mail: } psablonn@insa-rennes.fr

%----------------------------------------------------END------------------------------------------------------
\end{document}